# Designing Eco-Efficint Biowaste Valorisation Systems: Effects of Transport Distance and Facility Scale

Christian Ottini, Gwenola Yannou-Le Bris, Sandra Domenek and Felipe Buendia*

Université Paris-Saclay, INRAE, AgroParisTech, UMR SayFood, 91120 Palaiseau, France

felipe.buendia@inrae.fr

The increasing deployment of separate biowaste collection across Europe raises questions regarding environmentally and economically optimal treatment strategies. From this perspective, the balance between decentralised and centralised infrastructure is a potential lever to consider. The choice of technologies is another. Composting and anaerobic digestion (AD) offer viable circular bioeconomy pathways, but their performance depends on transport distance, facility size, seasonal waste fluctuations, and technological efficiency. Threshold distances and scale-dependent interactions affecting the environmental and economic performance of treatment plants remain poorly characterised. This study develops a unified techno-economic and life cycle assessment framework to quantify how transport constraints and treatment capacities shape biowaste valorisation systems. Three facility scales for composting and AD are evaluated across a range of transport distances and compared with local incineration in terms of environmental and economic performance. Transport burdens are modelled using vehicle life cycle assessment (LCA) and realistic heavy-duty truck data for electric-, gas-, and diesel-powered vehicles. Results reveal non-linear trade-offs between scale and transport, with zero-distance composting yielding the lowest combined impacts. Across all scenarios, 45 configurations outperform local incineration economically and environmentally, while 111 remain environmentally preferable despite higher costs, thus providing quantitative guidance for infrastructure planning.

## 1. Introduction

Biowaste constitutes the largest fraction of municipal solid waste and represents a key lever for achieving circular-economy targets through nutrient recovery, soil quality improvement, and renewable energy generation. Despite this potential, biowaste remains underutilised in many regions, with landfilling and incineration still dominating end-of-life pathways (Kaza et al., 2018). European regulatory developments, including the mandatory separate collection and landfill diversion requirements (European Parliament, Council of the European Union, 2018a, 2018b), are now accelerating the deployment of composting and anaerobic digestion (AD) infrastructures. Still, creating more sustainable biowaste systems depends heavily on local conditions and must take into account territorial, technological, and logistical interactions.

Among these interactions, transport plays a central role. Transport often accounts for a substantial share of operational costs and can represent a significant portion of environmental burdens, especially when biowaste is hauled over long distances (Lewerenz et al., 2023). Recent studies highlight how vehicle type, load factor, range autonomy, and operational conditions significantly influence transport emissions and costs (Sacchi et al., 2021; Syré and Göhlich, 2025). For biowaste management specifically, uncertainties in emission factors and logistical conditions further complicate the evaluation of transport impacts (Lewerenz et al., 2023). Consequently, determining whether biowaste should be treated locally or transported to larger valorisation hubs is a non-trivial decision.

Biological treatment technologies themselves exhibit scale-dependent behaviour and economic performance. Larger AD and composting facilities benefit from economies of scale (Tsilemou and Panagiotakopoulos, 2006), improved energy recovery options, and greater operational efficiencies. In contrast, small-scale units reduce transport requirements in highly urbanised areas where space constraints can be significant. Previous studies comparing centralised and decentralised AD configurations (Tian et al., 2021) demonstrate that the relative sustainability of each strategy depends on both technological performance and spatial organisation. However, most existing assessments are case-specific and rarely systematically explore the combined influence of transport distance, facility scale, and technology choice.

As a result, planning guidelines often rely on simplistic “maximum distance” thresholds that may overlook non-linear relationships between transport burdens and economies of scale. The absence of a generalisable, scenario-based understanding hampers municipalities and regional planners who must design infrastructures under evolving regulatory and territorial constraints.

To address these gaps, this study develops a unified techno-economic and environmental modelling framework that couples transport constraints with multi-scale facility design for biowaste valorisation. Instead of focusing on a single territory, the approach constructs a representative scenario space that combines three facility scales for composting and AD and evaluates them across a broad range of transport distances. Transport modelling is explicitly informed by recent advancements in heavy-duty vehicle LCA and realistic load-factor considerations (Sacchi et al., 2021), and operational costs (Basma and Rodríguez, 2023), while treatment performance is benchmarked against local incineration as a reference end-of-life option.
The objective of this study is to quantify how transport distance and facility scale jointly affect the environmental and economic performance of composting and anaerobic digestion systems, and to identify threshold distances at which decentralised or centralised configurations become preferable to local incineration.

## 2. Methodology

The following section outlines how transport modelling, facility performances, and scale-dependent factors are integrated into a unified framework to identify threshold distances and optimal configurations across the scenario space. A break-even point is defined as the value at which an alternative system becomes more eco-efficient than the base scenario.

### 2.1 Environmental modelling and system boundaries

The system environmental efficiencies were assessed using the LCA method. The impact assessment chosen was ReCiPe 2016 Endpoint (H) (Huijbregts et al., 2017). The normalisation factor used was World 2010 ("Normalization scores ReCiPe 2016 | RIVM," 2024), and the open source library Brightway (Mutel, 2017) for the assessment. The functional unit was defined as the transport and successive treatment of 1 tonne of biowaste. System boundaries extend from the point of collection (the bin) through transport to the facility and including all on-site treatment processes. Downstream application of compost or digestate was excluded. To ensure methodological consistency, the study adopts an attributional, cutoff-based LCA approach, as the life cycle inventories employed for modelling transport are of this type. The technology inventories were collected from our previous study (Ottini et al., 2025) which included the full treatment process, but no transport. Transportation was modelled using the Python library "*Carculator*" (Sacchi et al., 2021), with inventories migrated to Ecoinvent 3.11 for consistency. The chosen transport combinations were: diesel (ICEV-d), compressed gas (ICEV-g), and battery electric vehicle (BEV), for the 7.5t lorry and the year 2020, to keep consistency in the study.

### 2.2 Techno-Economic Analysis

The Techno-Economic Analysis (TEA) was developed using approximate cost functions from Tsilemou & Panagiotakopoulos (2006), who derived these equations using the capital expenditure (CAPEX) and operational expenditure (OPEX) costs of a range of established facilities distributed throughout Europe. The CAPEX values were updated to 2023 prices using the Chemical Engineering Plant Cost Index (CEPCI) ("CEPCI," 2014), while the OPEX values were inflated to 2023 using the Eurostat GDP implicit price deflator (Eurostat, 2025). This ensures that both capital and operational costs reflect current economic conditions across Europe. Moreover, the revenues for the sale of treatment-derived products, namely electricity, heat and compost, were accounted for each of the technologies producing one or more of these products. After expert consultation, the selling prices of the products were assumed to be 0.22 €/kWh for electricity, 0.075 €/kWh for heat, 50 €/t for the mature compost and a collection tax of 19 €/t. Transport costs were extracted from Basma & Rodríguez (2023) for the previously mentioned technologies in 2023.
The techno-economic performance is evaluated using the Net Present Value (NPV) formula, shown in Equation 1:

$$NPV_i = \sum_{t=1}^{T} \frac{REV_{i,t} - (OPEX_{i,t} + TRANS_{i,t})}{(1+r)^t} - CAPEX_i \quad (1)$$

where *i* denotes the technology under analysis, *t* the year of operation, *T* the project lifetime, and *r* the discount rate (assumed to be 3.5%).
$REV_{(i,t)}$ represents revenues from energy and material recovery, $OPEX_{(i,t)}$ the annual operating costs, and $TRANS_{(i,t)}$ the transport-related costs in year t. $CAPEX_i$ corresponds to the initial capital investment for technology *i*.

### 2.3 Scenario construction and assumptions

To explore the interplay between facility scale, technology, and transport distance, the study constructs a representative scenario space. Three facility scales are considered for both composting and anaerobic digestion (AD), reflecting typical small, medium, and large plants observed across Europe. Transport distances are varied systematically from 0 km (on-site treatment) to 150 km, capturing plausible ranges for local versus regional waste logistics. Seasonal variations in biowaste quantity are implicitly included via the annual average tonnage, and transport is modelled for diesel, biogas, and electric 7.5t trucks. Local incineration is included as a reference alternative. All scenarios assume full treatment capacity at each facility scale, and costs and environmental impacts are calculated per tonne of biowaste treated.

## 3. Results

The following section presents the techno-economic and environmental performance of the considered biowaste valorisation configurations. We first discuss the cost implications of centralised versus decentralised biotechnological alternatives, highlighting the economic break-even point relative to local incineration. Next, we assess the environmental impacts of the facilities, both with and without transport. Finally, configurations are compared in a Cost-Environmental impact space to identify an Eco-Efficient design zone relative to the proposed base scenario.

### 3.1 Cost Implications of a Centralised versus Decentralised System

Figure 1 presents the total treatment cost per tonne of biowaste for all evaluated configurations, combining three facility scales, three transport fuels, and transport distances. Large-scale incineration at 0 km is used as a reference and is shown as a horizontal threshold line.

At zero transport distance, all composting configurations and medium- and large-scale AD plants exhibit lower treatment costs than incineration, except for the small-scale AD facility, which remains less competitive due to higher operational expenditures. A clear economies-of-scale effect is visible across both technologies, with unit costs decreasing markedly from small to large installations.

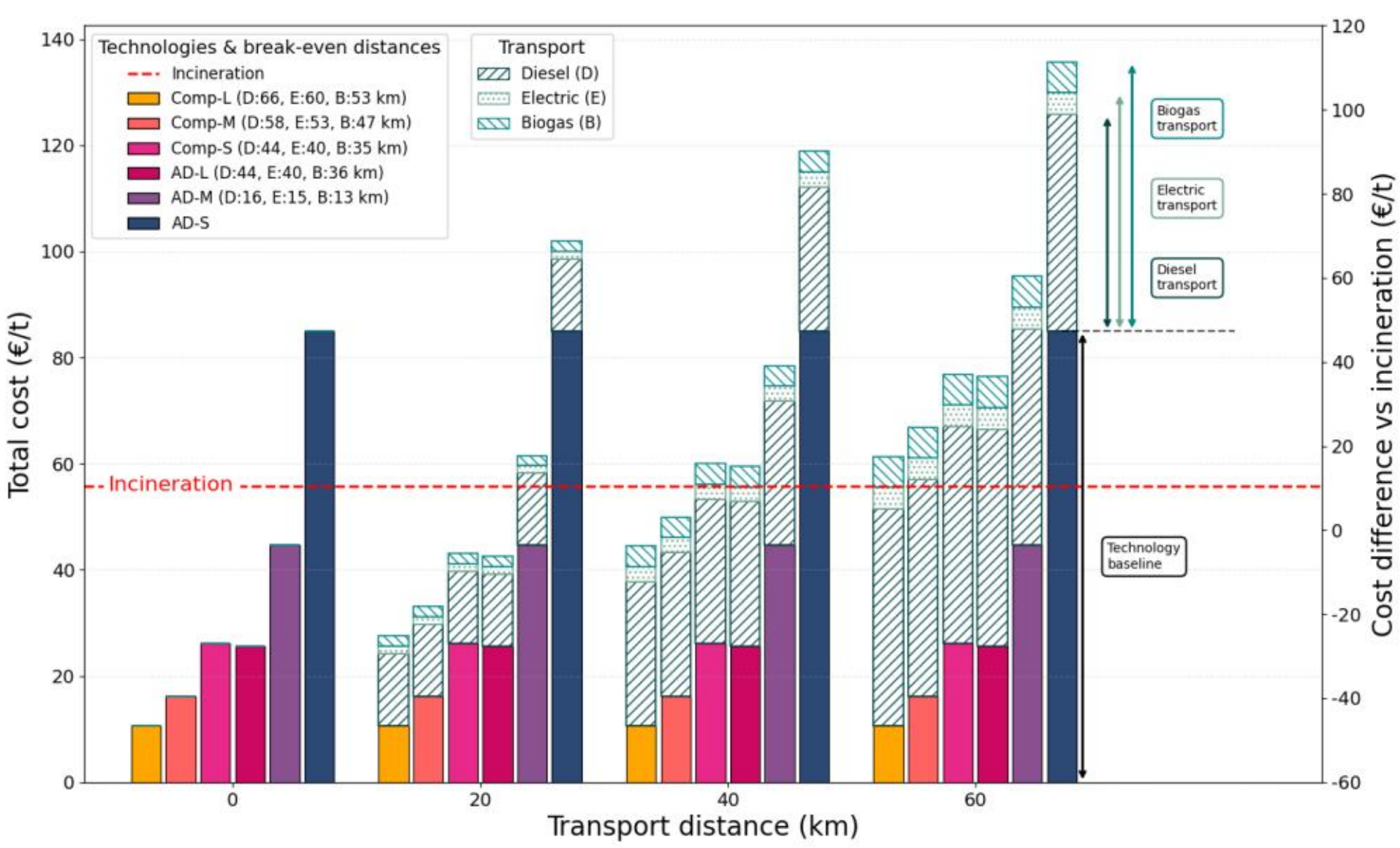

*Figure 1: Total treatment cost per tonne of biowaste for all technology-scale-transport configurations across transport distances. The red dashed line represents the cost of local large-scale incineration at 0 km. Comp: Composting, AD: Anaerobic Digestion, L, M, and S stands for large, medium, and small-sized facilities.*

As transport distance increases, costs rise linearly for all technologies and fuel types. Most configurations remain economically competitive up to approximately 20 km. Medium-scale AD reaches its break-even point slightly earlier, at approximately 16 km for diesel, 15 km for electric trucks, and 13 km for biogas trucks. Beyond 40 km, when considering diesel only, most technologies remain viable except those previously excluded. For alternative fuels, medium-scale AD and small composting facilities become less competitive than local incineration at 40 km for AD with electric trucks, 36 km with biogas, and at 40 km and 35 km for composting with electric and biogas trucks, respectively. At 60 km, only large composting with diesel remains competitive (up to 66 km), while medium composting is only competitive up to 58 km.

### 3.2 Environmental assessment of the configurations

Environmental results for all configurations are shown in Figure 2. Impacts are reported separately for treatment and the additional transport burden. Across all transport distances, alternative configurations are compared with the reference case of large-scale local incineration.

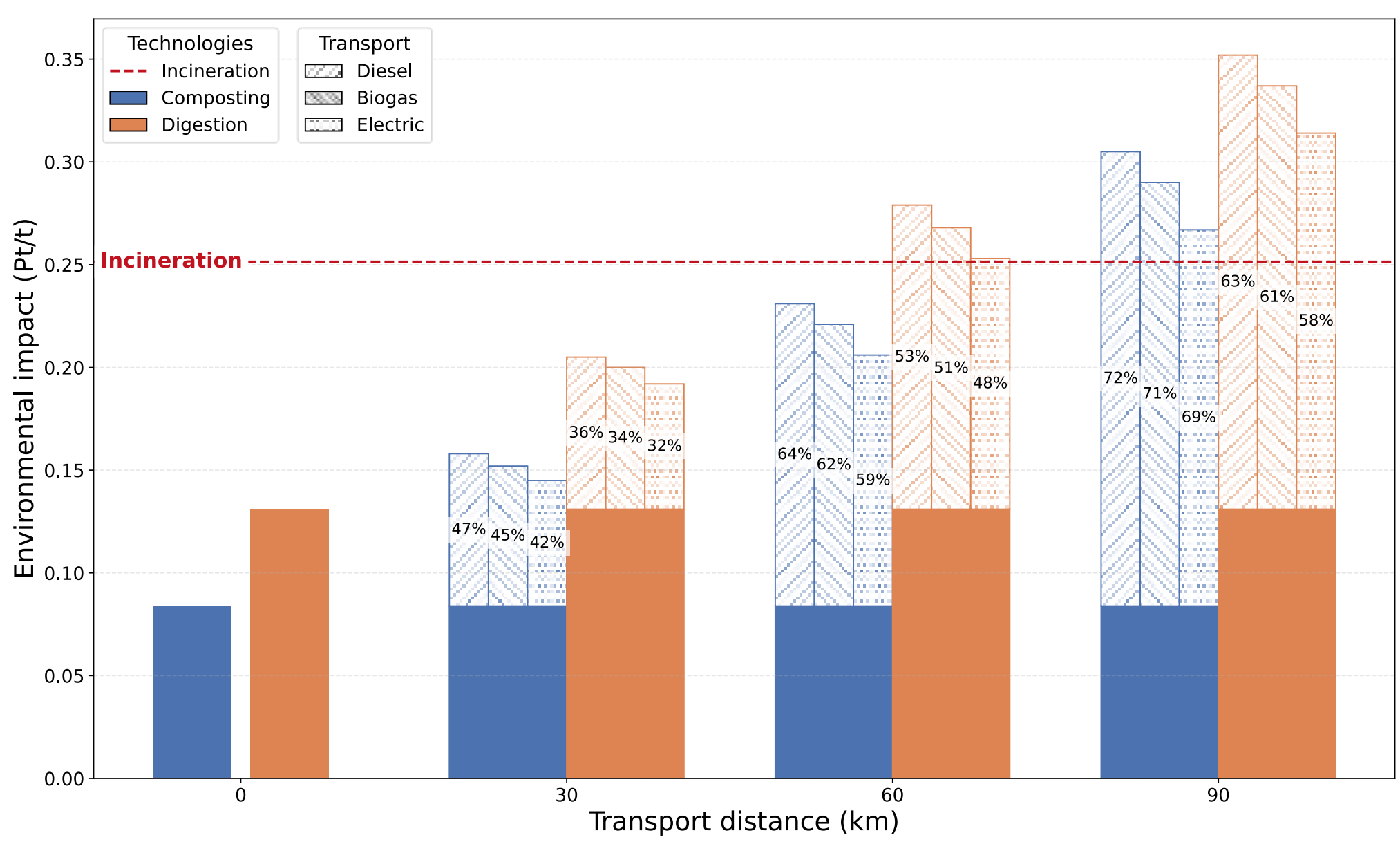

*Figure 2: Total environmental impact per tonne of biowaste for all technology-scale-transport configurations across transport distances. The red dashed line represents the impact of local large-scale incineration at 0 km. Solid bars show the impact of the treatment process, hatched bars indicate the transport-related impact, and the percentages in each hatched bar represent the share of transport in the total impact (treatment + transport).*

The results indicate that the environmental impact of alternative technologies is generally lower than that of incineration. Specifically, anaerobic digestion and composting have impacts approximately one-half and one-third that of incineration, respectively. At a transport distance of 30 km, all technologies still perform better environmentally than local incineration. However, at 60 km, all anaerobic digestion options become more environmentally burdensome than incineration. Finally, at 90 km, all alternative technologies result in higher impacts than local incineration of biowaste. The break-even points and the transport share corresponding to 50% of the total impact are reported in Table 1.

*Table 1: Configuration performance: break-even distances and transport distances at which transport accounts for 50% of the total environmental impact (treatment + transport).*

| Technology-transport mode | Break-even point (km) | Transport 50 % share (km) |
|---|---|---|
| Composting Diesel | 68.3 | 34 |
| Composting Biogas | 73.2 | 36.9 |
| Composting Electric | 82.3 | 41.2 |
| AD Diesel | 48.8 | 53.1 |
| AD Biogas | 52.7 | 57.4 |
| AD Electric | 59.2 | 64.4 |

### 3.3 Cost-Environmental Trade-Offs and Pareto-Optimal configurations

In the previous subchapters, we analysed the relationships and threshold results from the techno-economic and environmental assessments of the system. Here, we focus on integrating these two analyses using a graphical optimisation approach to identify configurations that achieve the best balance between economic and environmental performance. All scenarios were mapped in a cost versus environmental impact space. Figure 3 illustrates this coupling for selected transport distances and technologies. Configurations located in the lower-left quadrant, here called the “Eco-Efficient Design Zone”, represent solutions with both better environmental and economic performance compared to local incineration.

The results indicate that all composting systems reach Pareto efficiency at transport distances of less than 40 km for small facilities, 58 km for medium facilities, and up to 66 km for large facilities when diesel is the fuel of choice. With the same fuel choice, large AD facilities achieve optimality at intermediate transport distances (0–44 km), benefiting from economies of scale that offset increasing transport burdens. In contrast, medium-sized AD facilities are only competitive up to 16 km, and small AD facilities are never economically competitive and therefore do not appear in the optimal quadrant.

In the lower-right quadrant, we can still encounter sub-optimal solutions that are economically disadvantageous compared to local treatment through incineration, but remain environmentally valuable in terms of sole impact reduction.
Across the scenario space, large-scale composting with zero transport distance emerges as the most Eco-Efficient configuration (star in Fig. 3). This scenario achieves the lowest total treatment cost (10.65 € $t^{-1}$) and the lowest environmental impact (0.0839 Pt $t^{-1}$). Overall, 61 configurations outperform local incineration both economically and environmentally, while 111 remain environmentally preferable despite higher costs.

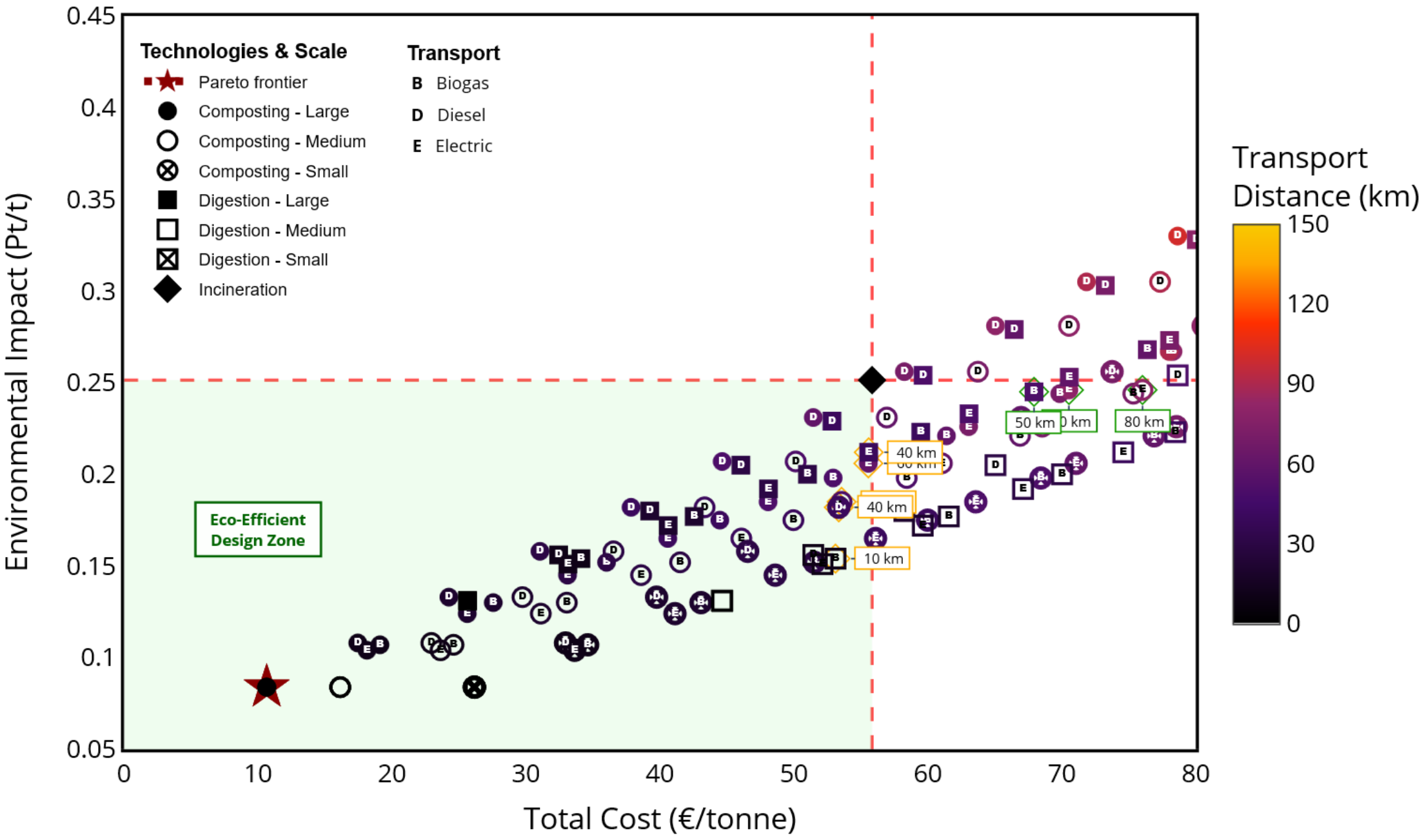

*Figure 3: Cost-Environmental Trade-Offs and Eco-Efficient design zone for biowaste treatment systems, including transportation integrated with anaerobic digestion and composting. The scatter plot shows the location of the various configurations in our research space.*

These results highlight non-linear relationships between facility scale, transport distance, and sustainability outcomes. Treating biowaste on-site or very close to the source (zero-distance treatment) consistently minimises combined environmental and economic impacts. This is primarily achieved by avoiding transport-related emissions and costs. This finding underscores the advantage of decentralised approaches in dense zones with reduced space, where many small facilities can efficiently manage local waste streams. Conversely, regional hubs may be justified where sufficient biowaste volumes allow large-scale facilities to exploit scale efficiencies, provided transport distances remain within identified thresholds.

## 4. Conclusions

This study shows that the environmental and economic performance of biowaste valorisation systems is driven by the interaction between facility scale, transport distance, and technology. Large-scale composting with zero transport distance is the most eco-efficient option, achieving the lowest cost (10.65 € $t^{-1}$) and environmental impact (0.0839 Pt $t^{-1}$). Small-scale AD is generally not competitive, while medium- and large-scale facilities benefit from economies of scale only up to specific distance thresholds. Alternative transport fuels further improve environmental performance and extend viable transport ranges.
The complexity and interdependence of technological, logistical, territorial, and social factors underscore the necessity of a holistic, systemic approach to biowaste management. Future analyses should integrate prospective LCA, techno-economic evaluations, and social considerations, including community acceptance, job creation, and local engagement. Systemic optimisation approaches are needed to identify non-dominated solutions that are simultaneously environmentally sustainable, economically viable, and socially equitable. Such forward-looking, systemic approaches are essential for designing robust, low-carbon, circular biowaste infrastructures that align with sustainable development goals in European cities and regions.
Overall, the proposed framework provides quantitative guidance for identifying when decentralised treatment is preferable and when regional hubs remain justified, supporting evidence-based planning of sustainable biowaste infrastructures.